\documentclass[reqno]{amsart}
\usepackage{amsmath,amssymb,amsthm,ascmac}
\usepackage[mathscr,mathcal]{eucal}
\usepackage{latexsym}
\usepackage{enumerate}
\usepackage[dvipdfmx]{graphicx}
\usepackage{lineno,amscd,epsfig,color,ulem}
\usepackage{pb-diagram} 

\begin{document}




\theoremstyle{definition}
\newtheorem{defi}{\textbf{Definition}}[section]
\newtheorem{thm}[defi]{\textbf{Theorem}}
\newtheorem{lem}[defi]{\textbf{Lemma}}
\newtheorem{prop}[defi]{\textbf{Proposition}}
\newtheorem{cor}[defi]{\textbf{Corollary}}
\newtheorem{ex}[defi]{\textbf{Example}}
\newtheorem{rem}[defi]{\textbf{Remark}}
\newtheorem*{corr}{\textbf{Corollary}}

\theoremstyle{plain}
\newtheorem{maintheorem}{Theorem}
\newtheorem{theorem}{Theorem }[section]
\newtheorem{proposition}[theorem]{Proposition}
\newtheorem{mainproposition}{Proposition}
\newtheorem{lemma}[theorem]{Lemma}
\newtheorem{corollary}[theorem]{Corollary}
\newtheorem{maincorollary}{Corollary}
\newtheorem{claim}{Claim}
\renewcommand{\themaintheorem}{\Alph{maintheorem}}
\theoremstyle{definition} \theoremstyle{remark}
\newtheorem{remark}[theorem]{Remark}
\newtheorem{example}[theorem]{Example}
\newtheorem{definition}[theorem]{Definition}
\newtheorem{problem}{Problem}
\newtheorem{question}{Question}
\newtheorem{exercise}{Exercise}

\newtheorem*{subject}{\UTF{FFFD}\UTF{0693}I}
\newtheorem*{mondai}{Problem}
\newtheorem{lastpf}{\CID{1466}\UTF{FFFD}}
\newtheorem{lastpf1}{proof of theorem2.4}
\renewcommand{\thelastpf}{}
\renewcommand{\labelenumi}{(\Roman{enumi})}
\newcommand{\Add}[1]{\textcolor{red}{#1}}
\newcommand{\Erase}[1]{\textcolor{red}{\sout{{#1}}}}
\newcommand{\ul}{\underline}

\newtheorem{last}{Theorem}

\renewcommand{\thelast}{}
\renewcommand{\proofname}{\textup{Proof.}}

\renewcommand{\theequation}{\arabic{section}.\arabic{equation}}
\makeatletter
\@addtoreset{equation}{section}

\title[]
{Density of periodic measures and large deviation principle for generalized ($\alpha,\beta$)- transformations}


\author[M. Shinoda]{Mao Shinoda}
\address{Department of Mathematics \\
Ochanomizu University \\ 
2-1-1, Otsuka, Bunkyo-ku,Tokyo 112-8610, JAPAN}
\email{shinoda.mao@ocha.ac.jp}
\author[K. Yamamoto]{Kenichiro Yamamoto}
\address{Department of General Education \\
Nagaoka University of Technology \\
Nagaoka 940-2188, JAPAN
}
\email{k\_yamamoto@vos.nagaokaut.ac.jp}


\date{}

\maketitle
\large

\begin{abstract}

We introduce generalized $(\alpha,\beta)$-transformations,
which include all $(\alpha,\beta)$
and generalized $\beta$-transformations,
and prove that all transitive
generalized $(\alpha,\beta)$-transformations
satisfy the level-2 large deviation principle
with a unique measure of maximal entropy.
A crucial step in our proof is to establish density of periodic measures in the set of ergodic measures.


\end{abstract}

\section{Introduction}

In this paper, we consider piecewise monotonic maps
on the unit interval $[0,1]$.
We say that $T\colon [0,1]\to [0,1]$
is a \textit{piecewise monotonic map} if
there exist integer $k>1$ and
$0=c_0<c_1<\cdots<c_k=1$, which we call the
\textit{critical points},
such that $T|_{(c_{i-1}, c_i)}$ is strictly monotonic and continuous for each $1\leq i\leq k$.
Throughout this paper, we further assume the
following conditions for a piecewise monotonic map $T$.
\begin{itemize}
\item 
The \textit{topological entropy} $h_{\rm top}(T)$ of $T$
is positive
(see \cite[Ch. 9]{BB} for the definition of topological
entropy for piecewise monotonic maps).
\item
$T$ is \textit{transitive}, i.e.,
there exists a point $x\in [0,1]$ whose
forward orbit $\{T^n(x):n\ge 0\}$ is dense
in $[0,1]$.
\end{itemize}
Under these conditions, it is proved in \cite[Theorem 4]{H2} that 
there exists a unique measure of maximal entropy for $T$,
that is, a $T$-invariant measure whose metric entropy coincides with $h_{\rm top}(T)$.
The aim of this paper is to investigate
whether the large deviation principle holds for
a piecewise monotonic map with
the unique measure of maximal entropy as a reference.

Let $\mathcal{M}([0,1])$ be the set of all
Borel probability measures on $[0,1]$ endowed with
the weak${}^{\ast}$-topology.
We say that $([0,1],T)$ satisfies
the \textit{(level-2) large deviation principle}
with the unique measure of maximal entropy $m$
as a reference
if there exists a lower semi-continuous function $\mathcal{J}\colon \mathcal{M}([0,1])\to [0,\infty]$, called
a \textit{rate function}, such that
\[\limsup_{n\rightarrow\infty}\frac{1}{n}\log m\left(\left\{x\in [0,1]:\frac{1}{n}\sum_{j=0}^{n-1}\delta_{T^j(x)}\in\mathcal{K}\right\}\right)\le
-\inf_{\mathcal{K}}\mathcal{J} \] 
holds for any closed set
$\mathcal{K}\subset\mathcal{M}([0,1])$ and
\[\liminf_{n\rightarrow\infty}\frac{1}{n}\log m\left(\left\{x\in [0,1]:\frac{1}{n}\sum_{j=0}^{n-1}\delta_{T^j(x)}\in\mathcal{U}\right\}\right)\ge
-\inf_{\mathcal{U}}\mathcal{J}\]
holds for any open set $\mathcal{U}\subset\mathcal{M}([0,1])$.
Here $\delta_y$ signifies the Dirac mass
at a point $y\in [0,1]$.
We refer to \cite{E} for a general theory of large deviations and its background
in statistical mechanics.

In this literature, the case that $T$ has a constant slope is important because this condition
implies that the unique measure of maximal entropy
is absolutely continuous to the Lebesgue measure.
Hence we focus our attention to piecewise monotonic maps with a constant slope $\beta>1$.
A familiar example included in this family is a $\beta$-transformation.
\begin{itemize}
    \item 
    {\bf $\beta$-transformations.}\\
The $\beta$-transformation $T_{\beta}\colon [0,1]\to [0,1]$ with $\beta>1$ was introduced by R\'enyi \cite{Re} and defined by
\[T_{\beta}(x)=
\left\{
\begin{array}{ll}
\beta x\ (\text{mod}\ 1) & (x\not=1),\\
\displaystyle\lim_{y\rightarrow 1-0}(\beta y\ (\text{mod}\ 1)) &  (x=1).
\end{array}
\right.
\]
\end{itemize}
This family has recently attracted attention in the setting  beyond specification, 
since $T_{\beta}$ does not satisfy the
specification property for Lebesgue almost
parameter $\beta>1$ (\cite{Bu2,S}).
Inspired by $\beta$-transformations, many authors have considered various generalizations of $T_{\beta}$:
{\it $(\alpha,\beta)$-transformations}
(\cite{CLR,CY,FP,H5,RS}), {\it $(-\beta)$-transformations} (\cite{IS,LS,N,SY}), {\it generalized $\beta$-transformations} (\cite{FP,G,Su,Th}).
\begin{itemize}
    \item 
    {\bf $(\alpha,\beta)$-transformations.}\\
The $(\alpha,\beta)$-transformation $T_{\alpha,\beta}\colon [0,1]\to [0,1]$ with $\beta>1$ and
$0\le \alpha<1$ was introduced by Parry (\cite{P}) and defined by
\[
T_{\alpha,\beta}(x)=
\left\{
\begin{array}{ll}
\beta x+\alpha\ (\text{mod}\ 1) & (x\not=1),\\
\displaystyle\lim_{y\rightarrow 1-0}(\beta y+\alpha\ (\text{mod}\ 1)) &  (x=1).
\end{array}
\right.
\]
    
    \item
    {\bf Generalized $\beta$-transformations.}\\
    Let $\beta>1$ and $k$ be a smallest integer
    not less than $\beta$,
    fix $E=(E_1,\ldots,E_k)\in\{+1,-1\}^k$,
    and consider $k$ intervals
    \[I_1:=[0,1/\beta),\ I_2:=[1/\beta,2/\beta),\ldots,
    I_k:=[k-1/\beta,1].\]
    The generalized $\beta$-transformation
    $T_{\beta,E}\colon [0,1]\to [0,1]$ was
    introduced by G\'ora (\cite{G}) and defined by
    \[
T_{\beta,E}(x)=
\left\{
\begin{array}{ll}
\beta x-i+1 & (x\in I_i,\ E(i)=+1),\\
-\beta x+i & (x\in I_i,\ E(i)=-1).
\end{array}
\right.
\]
If $E=(-1,\ldots,-1)$, then
we call $T_{\beta,E}$ a $(-\beta)$-transformation.
\end{itemize}

Pfister and Sullivan (\cite{PS}) established the large deviation principle for $\beta$-transformations for any $\beta>1$, 
which is the first work on the large deviation principle in this family without the specification property.
In \cite{CY}, Chung and the second author
proved that the large deviation principle holds for
$(\alpha,\beta)$ and generalized $\beta$-transformations
in the following parameters:
\begin{itemize}
    \item 
    $(\alpha,\beta)$-transformations for
    $0\le\alpha<1$ and $\beta>2$.
    
    \item
    Generalized $\beta$-transformations for
    $\frac{1+\sqrt{5}}{2}<\beta<2$ and $E=(-1,-1)$.
\end{itemize}
Then it is natural to ask whether does the large deviation principle hold for the other parameters.
In this paper not only we get an affirmative answer to the all parameters, but also we show that it holds for more general class of piecewise monotonic maps with
constant slope, which we call
generalized $(\alpha,\beta)$-transformations.
\begin{itemize}
    \item
    {\bf Generalized $(\alpha,\beta)$-transformations.}
Let $0\le \alpha<1$, $\beta>1$,
$k$ be a smallest integer not less than
$\alpha+\beta$,
fix $E=(E_1,\ldots,E_k)\in\{+1,-1\}^k$, 
and consider $k$ intervals
    \[\ \ \ \ \ \ \ \ I_1:=\left[0,\frac{1-\alpha}{\beta}\right),\ I_2:=\left[\frac{1-\alpha}{\beta},\frac{2-\alpha}{\beta}\right),\ldots,
    I_k:=\left[\frac{k-1-\alpha}{\beta},1\right].\]
We define the generalized $(\alpha,\beta)$-transformation
$T_{\alpha,\beta,E}\colon [0,1]\to [0,1]$ by
\[
T_{\alpha,\beta,E}(x)=
\left\{
\begin{array}{ll}
\alpha+\beta x-i+1 & (x\in I_i,\ E(i)=+1),\\
-\alpha-\beta x+i & (x\in I_i,\ E(i)=-1).
\end{array}
\right.
\]
\end{itemize}
Generalized $(\alpha, \beta)$-transformations are clearly a generalization of both
$(\alpha, \beta)$-transformations and generalized $\beta$-transformations.
The graphs of $T_{\beta}$, $T_{\alpha,\beta}$,
$T_{\beta,E}$ and $T_{\alpha,\beta,E}$
are plotted in Figure 1.

Now we state our first main result of this paper. 

\begin{figure}
    \centering
    \includegraphics[width=300pt]{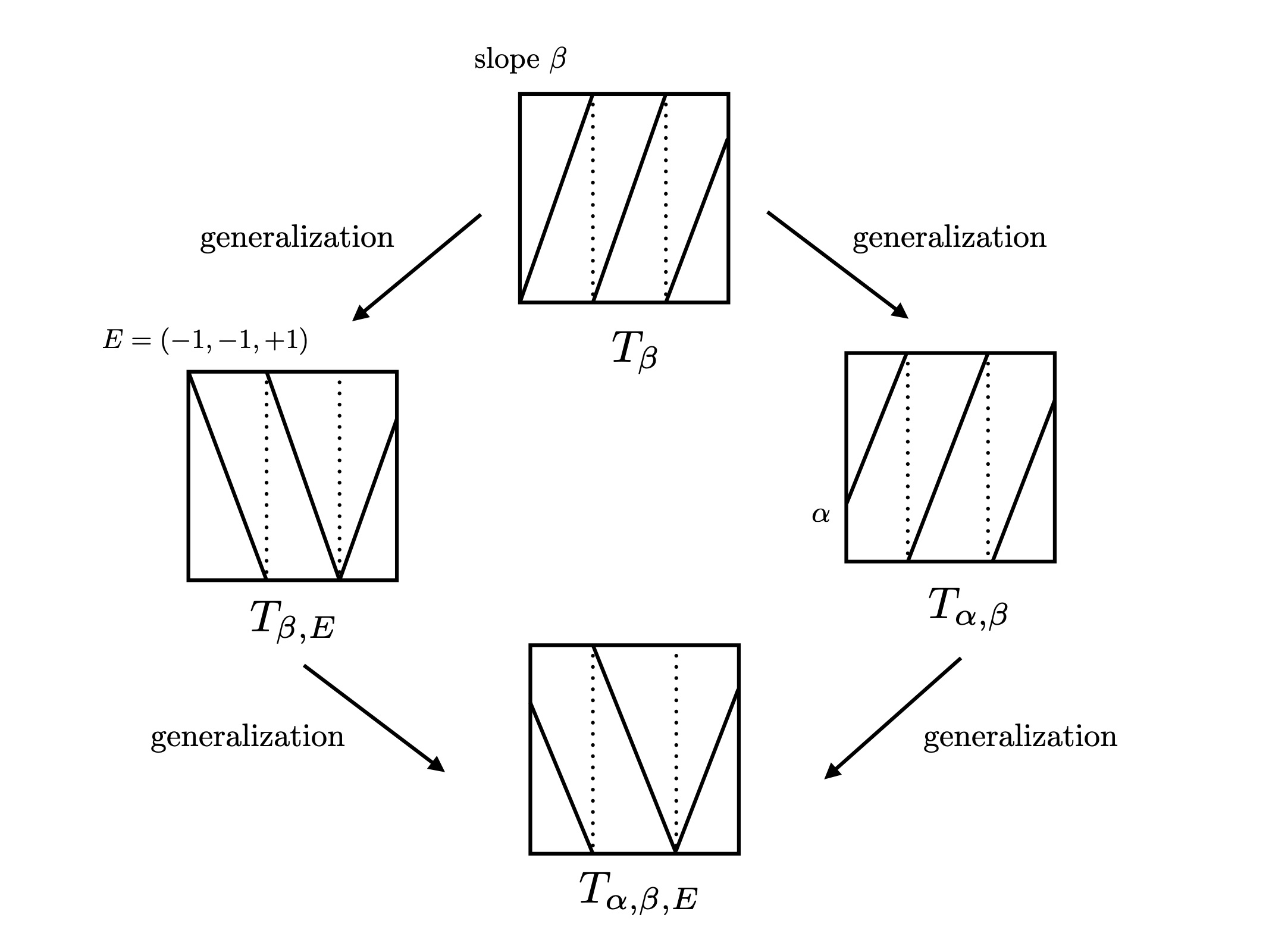}
    \hspace{-1cm}\caption{Graphs of $T_{\beta}$, $T_{\alpha,\beta}$,
    $T_{\beta,E}$ and $T_{\alpha,\beta,E}$.}
    \label{generalization}
\end{figure}

\begin{maintheorem}
\label{main1}
Let $T\colon [0,1]\to [0,1]$ be a
transitive generalized $(\alpha,\beta)$-\\
transformation.
Then $([0,1],T)$ satisfies the level-2
large deviation principle with the unique measure of
maximal entropy.
\end{maintheorem}



A crucial step to prove Theorem \ref{main1} is to show  
density of periodic measures.
For a metrizable space $X$ and
a Borel measurable map $f\colon X\to X$,
denote by $\mathcal{M}(X)$ the set of all
Borel probability measures on $X$ endowed with the
weak${}^{\ast}$-topology, by
$\mathcal{M}_f(X)\subset \mathcal{M}(X)$
the set of $f$-invariant ones, and
by $\mathcal{M}^e_f(X)\subset\mathcal{M}_f(X)$
the set of ergodic ones.
We say that $\mu\in\mathcal{M}(X)$ is a
\textit{periodic measure} if there exist $x\in X$ and $n>0$ such that
$f^n(x)=x$ and $\mu=\delta_n^f(x):=1/n\sum_{j=0}^{n-1}\delta_{f^j(x)}$ hold. Then, it is clear that $\mu\in\mathcal{M}_f^e(X)$.
We denote by $\mathcal{M}_f^p(X)\subset\mathcal{M}_f^e(X)$
the set of all periodic measures on $X$.
It is established by Chung and the second author that
 the level-2 large deviation principle for a piecewise monotonic map 
 is followed by
 density of periodic measures in the set of
 ergodic measures
 with the irreducibility of its Markov diagram 
 (see \cite[Theorem A]{CY}).
Although the irreducibility
is slightly stronger than the transitivity,
we show in this paper that this result
remains true if the irreducibility of a Markov diagram is replaced by the transitivity of a map.


\vspace{0.2cm}

\noindent
{\bf Proposition A.}
\textit{Let $T\colon [0,1]\to [0,1]$ be
a transitive piecewise monotonic map with
$h_{{\rm top}}(T)>0$.
Suppose that $\mathcal{M}^p_T([0,1])$ is dense in $\mathcal{M}_T^e([0,1])$.
Then $([0,1],T)$ satisfies the level-2
large deviation principle with the unique measure of
maximal entropy.}
\vspace{0.1cm}

Hence Theorem \ref{main1} follows from
Proposition A and the following theorem,
which is the second main result of this paper:
\begin{maintheorem}
\label{diagram}
Let $T\colon [0,1]\to [0,1]$ be a
transitive generalized
$(\alpha,\beta)$-\\
transformation.
Then $\mathcal{M}_T^p([0,1])$ is dense in
$\mathcal{M}_T^e([0,1])$.
\end{maintheorem}
Our proof of Theorem B is based on the work by Hofbauer and Raith (\cite{HR}) where
density of periodic measures was proved for a piecewise monotonic map consisting of two monotonic pieces. 
Since the family of generalized $(\alpha, \beta)$-transformations includes maps with more than three monotonic pieces, we need to improve their method. 
One key difference between \cite{HR} and ours is Proposition \ref{no long edges}, which is one of the novelty of this paper (see also Remark \ref{novelty}).


The remainder of this paper is organized as follows.
In \S2, we establish our definitions and
prepare several facts.
Subsequently, we present proofs of
Proposition A and Theorem \ref{diagram} in \S3.


\section{Preliminaries}
\subsection{Symbolic dynamics}
Let $\mathbb{N}_0$ be the set of non-negative integers.
For a finite or countable set $A$,
we denote by $A^{\mathbb{N}_0}$ the one-sided
infinite product of $A$ equipped with the product topology of the discrete topology of $A$.
To simplify the notation,
given integers
$i\le j$ and $x_i,x_{i+1},\ldots,x_j\in A$,
we set $x_{[i,j]}:=x_ix_{i+1}\cdots x_j$.
We also write
$x_{[i,j)}:=x_{[i,j-1]}$ and
similarly for $x_{(i,j]}$ and $x_{(i,j)}$.
A sequence $\ul{x}\in A^{\mathbb{N}_0}$ will be defined by
all of its coordinates $x_n\in A$ with $n\ge 0$.
Let $\sigma$ be the shift map on $A^{\mathbb{N}_0}$ 
 (i.e.,~$(\sigma (\ul{x}))_n= x_{n+1}$ for each $n\ge 0$ and $\ul{x}\in A^{\mathbb{N}_0}$). 
When a subset  $\Sigma^+$ of $A^{\mathbb{N}_0}$ is $\sigma$-invariant and  closed,  
we call it 
a \textit{subshift}
and call $A$ the \textit{alphabet} of 
$\Sigma^+$.
For a matrix $M= (M_{ij})_{(i,j)\in A^2}$, each entry of which is $0$ or $1$, we define a subshift
$\Sigma^+_M\subset A^{\mathbb{N}_0}$ by
\[
\Sigma^+_M =\{ \ul{x}\in A^{\mathbb{N}_0} : \text{$M_{x_n x_{n+1}} =1 $ for all $n\ge 0$}\}
\]
and call $\Sigma^+_M$ a \textit{Markov shift}
with an \textit{adjacency matrix} $M$.

For a subshift $\Sigma^+$ on an alphabet $A$,
we set $\mathcal{L}(\Sigma^+):=\{x_{[0,n]}:\ul{x}\in\Sigma^+,\ n\ge 0\}$ and
$[u] := \left\{ \ul{x}\in \Sigma^+ : u=x_{[0,|u|)}\right\}$ for each $u\in\mathcal{L}(\Sigma^+)$, where $|u|$
denotes the length of $u$.
A word $v\in\mathcal{L}(\Sigma^+)$ is called a subword
of $u=u_0\cdots u_n\in\mathcal{L}(\Sigma^+)$ if
$v=u_{[i,j]}$ for some $0\le i\le j\le n$.
For $u,v\in\mathcal{L}(\Sigma^+)$, we use
juxtaposition $uv$ to denote the word
obtained by
the concatenation and $u^{\infty}$ means a one-sided infinite sequence
$uuu\cdots\in A^{\mathbb{N}_0}$.
Moreover, we set
$\mathcal{S}(u):=\{i\in A:ui\in\mathcal{L}(\Sigma^+)\}$ for
$u\in\mathcal{L}(\Sigma^+)$.
Finally, we say that $\Sigma^+$ is \textit{transitive} if for any $u,v\in\mathcal{L}(\Sigma^+)$,
we can find $w\in\mathcal{L}(\Sigma^+)$ such that $uwv\in\mathcal{L}(\Sigma^+)$ holds.
For the rest of this paper, we denote by
\textcolor{black}{$h_{\sigma}(\mu)$} the \textit{metric entropy} of $\mu\in\mathcal{M}_{\sigma}(\Sigma^+)$.

\subsection{Markov diagram}

Let \textcolor{black}{$X =[0,1]$ and} $T\colon X\to X$ be a
piecewise monotonic map with
\textcolor{black}{critical points $0=c_0<c_1<\cdots<c_k=1$.}
\textcolor{black}{Let $X_T:=\bigcap_{n=0}^{\infty}T^{-n}(\bigcup_{j=1}^k (c_{j-1},c_j))$, and}
define the \textit{coding map} $\mathcal{I}\colon X_T\to \{1,\ldots,k\}^{\mathbb{N}_0}$ by
\[(\mathcal{I}(x))_n=j\text{ if and only if }T^n(x)\in (c_{j-1},c_j),\]
which is injective since $T$ is transitive
(see \cite[Proposition 6.1]{Y2}).
We denote the closure of $\mathcal{I}(X_T)$ \textcolor{black}{in $\{1,\ldots,k\}^{\mathbb{N}}$} by $\Sigma_T^+$.
Then, $\Sigma_T^+$ is a subshift,
and $(\Sigma_T^+,\sigma)$ is called
the \textit{coding space} of $(X,T)$.
We use the following notations:
    
\begin{itemize}
\item 
$\displaystyle\ul{a}^{(i)}:=\lim_{x\rightarrow c_i+0}
\mathcal{I}(x)$,
$\displaystyle\ul{b}^{(i)}:=\lim_{x\rightarrow c_i-0}
\mathcal{I}(x)$  for $1\le i\le k-1$.
\item
${\rm adj}(\ul{a}^{(i)}):=\ul{b}^{(i)}$,
${\rm adj}(\ul{b}^{(i)}):=\ul{a}^{(i)}$
for $1\le i\le k-1$.
\item
$\displaystyle\ul{a}:=\lim_{x\rightarrow +0}
\mathcal{I}(x)$,
$\displaystyle\ul{b}:=\lim_{x\rightarrow 1-0}
\mathcal{I}(x)$.
\end{itemize}
We also set $\mathcal{CR}:=\{\ul{a}^{(i)},\ul{b}^{(i)}:
1\le i\le k-1\}$ and call it a \textit{critical set}.

In what follows, we define the Markov diagram, introduced by Hofbauer (\cite{H2}), which is a countable oriented graph
with subsets of $\Sigma_T^+$ as vertices.
Let $C\subset\Sigma_T^+$ be a closed subset with $C\subset [j]$ for some
\textcolor{black}{$1\le j\le k$}.
We say that a non-empty closed subset $D\subset\Sigma_T^+$ is a \textit{successor} of $C$ if
$D=[l]\cap\sigma(C)$ for some
\textcolor{black}{$1\le l\le k$}.
The expression $C\rightarrow D$
denotes that
$D$ is a successor of $C$.
Now, we define a set $\mathcal{D}_T$ of vertices by induction. First, we set
$\mathcal{D}_0:=\{[1],\ldots,[k]\}$. If $\mathcal{D}_n$ is defined for $n\ge 0$, then we set
\[\mathcal{D}_{n+1}:=\mathcal{D}_n\cup \{D:\textcolor{black}{D\text{ is a successor for some }C\in\mathcal{D}_n}\}.\]
We note that $\mathcal{D}_n$ is a finite set for each $n\ge 0$ since the number of successors of any
closed subset of $\Sigma_T^+$ is
at most $k$ by definition. Finally, we set
\[\mathcal{D}_T:=\bigcup_{n\ge 0}\mathcal{D}_n.\]
The oriented graph $(\mathcal{D}_T,\rightarrow)$ is called
the \textit{Markov diagram} of $T$.
For notational simplicity, we use the expression $\mathcal{D}$ instead of $\mathcal{D}_T$
if no confusion arises.

For $\ul{x}\in\Sigma_T^+$ and $n\ge 0$, we set
$D_n^{\ul{x}}:=\sigma^n(x_{[0,n]})$.
We define a sequence $\{R_m^{\ul{x}}\}_{m\ge 0}$ of
integers inductively as follows.
First, we set $R_0^{\ul{x}}:=0$.
If $R_m^{\ul{x}}$ is defined for $m\ge 0$,
then let
\[R_{m+1}^{\ul{x}}:=\min\{n>R_m:\#(\mathcal{S}(x_{[0,n-1)}))\ge 2\}.\]
We also set $r_m^{\ul{x}}:=R_m^{\ul{x}}-R_{m-1}^{\ul{x}}$.
Now, we summarize properties of
Hofbauer's Markov Diagram, which are
appeared in \cite{H3} (see also \cite[Page 224]{HR}).

\begin{prop}
\label{diagram-3}
Let $(\mathcal{D},\rightarrow)$ be the Markov diagram
of $T$.
\vspace{0.2cm}\\
(1) $\mathcal{D}=\{D_n^{\ul{x}}:\ul{x}\in\mathcal{CR}
\cup\{\ul{a},\ul{b}\},\ n\ge 0\}$.
\vspace{0.2cm}\\
(2) $(\mathcal{D},\rightarrow)$ has the following arrows:
\begin{itemize}
\item
For any $\ul{x}\in\mathcal{CR}\cup\{\ul{a}, \ul{b}\}$ and any $n\ge 0$,
$D_n^{\ul{x}}\rightarrow D_{n+1}^{\ul{x}}$.

\item
For any $\ul{x}\in\mathcal{CR}\cup\{\ul{a}, \ul{b}\}$ and any $m\ge 1$,
$D_{R_m^{\ul{x}}-1}^{\ul{x}}\rightarrow D_{r_m^{\ul{x}}}^{f_m(\ul{x})}$.
Here $f_m(\ul{x})$ is a unique point in $\mathcal{CR}$
so that $f_m(\ul{x})\not=\sigma^{R_{m-1}^{\ul{x}}}(\ul{x})$
and $f_m(\ul{x})\in D_{R_{m-1}}^{\ul{x}}$ hold.
\end{itemize}
(3) Let $\ul{x}\in\mathcal{CR}\cup\{\ul{a}, \ul{b}\}$.
\begin{itemize}
    \item 
    For any $0\le n\le r_m^{\ul{x}}-1$, we have
$D^{\ul{x}}_{n+R_{m-1}}\subset D_n^{f_m(\ul{x})}$.
In particular,
$x_{[R_{m-1}^{\ul{x}},R_m^{\ul{x}})}=f_m(\ul{x})_{[0,r_m^{\ul{x}})}$ holds.
    \item
    There exists an integer $q$ such that
    $r_m^{\ul{x}}=R_q^{f_m(\ul{x})}$.
\end{itemize}
(4) If $C\in\mathcal{D}$ and $\#\mathcal{S}(C)>2$, then $\mathcal{S}(C)\cap\mathcal{D}_0\not=\emptyset$.
\end{prop}

For a subset $\mathcal{C}\subset \mathcal{D}$, we define a matrix
$M(\mathcal{C})=(M(\mathcal{C})_{C,D})_{(C,D)\in\mathcal{C}^2}$ by
\[ M(\mathcal{C})_{C,D}=
\left\{
\begin{array}{ll}
1 & (\text{$C\rightarrow D$}), \\
0 & (\text{otherwise}).
\end{array}
\right. \]
Then, $\Sigma^+_{M(\mathcal{C})}=\{\ul{C}\in\mathcal{C}^{\mathbb{N}_0}:C_n\rightarrow
C_{n+1},n\ge 0\}$ is a \textcolor{black}{one-sided}
Markov shift with a countable alphabet $\mathcal{C}$
and an adjacency matrix $M(\mathcal{C})$.
For notational simplicity,
we denote $\Sigma_{\mathcal{C}}^+$ instead of $\Sigma^+_{M(\mathcal{C})}$.
We say that $\mathcal{C}$
is \textit{irreducible} if for any $C,D\in\mathcal{C}$,
there are finite vertices $C_0,\ldots,C_n\in\mathcal{C}$
such that $C_i\rightarrow C_{i+1}$ for $0\le i\le n-1$
(i.e.,
$C_{[0,n]}\in\mathcal{L}(\Sigma^+_{\mathcal{C}})$),
$C_0=C$ and $C_n=D$,
and if every subset of $\mathcal{D}$,
which contains
$\mathcal{C}$ does not have this property.
It is clear that $\Sigma^+_{\mathcal{C}}$
is transitive if and only if
$\mathcal{C}$ is
irreducible.
We define a map
$\Psi\colon \Sigma_{\mathcal{D}}^+\to \{1,\ldots,k\}^{\mathbb{N}_0}$ by
\[ \Psi((C_n)_{n\in\mathbb{N}_0}):=(x_n)_{n\in\mathbb{N}_0}\text{ for }(C_n)_{n\in\mathbb{N}_0}\in\Sigma^+_{\mathcal{D}},\]
where
\textcolor{black}{$1\le x_n\le k$} is a unique integer such that $C_n\subset [x_n]$ holds for each $n\in\mathbb{N}_0$.
Then it is not difficult to see that
$\Psi$ is continuous, countable to one and
$\Psi(\Sigma_{\mathcal{D}}^+)=\Sigma_T^+$.
We say that $\mu\in\mathcal{M}^e_{\sigma}(\Sigma_T^+)$
is \textit{liftable} if
there is $\overline{\mu}\in\mathcal{M}^e_{\sigma}
(\Sigma_{\mathcal{D}}^+)$
such that $\mu=\overline{\mu}\circ\Psi^{-1}$ holds.
It is known that not every
$\mu\in\mathcal{M}^e_{\sigma}(\Sigma_T^+)$ is liftable
in general although $\Psi$ is surjective (see \cite{K}).
In \cite{H2}, Hofbauer provided a sufficient
condition for the liftability.

\begin{lem}(\cite[Lemma 3]{H2})
\label{positive}
If $\mu\in\mathcal{M}_{\sigma}^e(\Sigma_T^+)$ has positive metric entropy,
then $\mu$ is liftable.
\end{lem}

We recall two important facts for
transitive piecewise monotonic maps.

\begin{lem}(\cite[Theorem 11]{H3},\cite[Page 224]{HR})
\label{tran}
There exists an irreducible subset $\mathcal{C}\subset\mathcal{D}$ satisfying the
following properties:
\begin{itemize}
\item
$\Psi(\Sigma_{\mathcal{C}}^+)=\Sigma_T^+$.

\item
$C\in\mathcal{C}$ and $C\rightarrow D$ implies
$D\in\mathcal{C}$.

\item
There exists an integer $n_0$ such that
$D_{n_0}^{\ul{x}}\in\mathcal{C}$ holds
for any $\ul{x}\in\mathcal{CR}\cup\{\ul{a}, \ul{b}\}$.
\end{itemize}
\end{lem}

\begin{thm}(\cite[Theorem 1]{HR})
\label{HR-thm}
Suppose that
there are integers $N_0$ and $N_1$ such that
for any $\ul{x}\in\mathcal{CR}$,
and any $j\in\mathbb{N}$ with $R_j^{\ul{x}}>N_0$,
there exist an integer $1\le m<j$,
a periodic point $\ul{p}$ with period
$l$ and $u\in\mathcal{L}(\Sigma_T^+)$
with
$|u|\ge R_j^{\ul{x}}-R_m^{\ul{x}}-N_1$
such that
$u$ is a subword of both $x_{[R_m^{\ul{x}},R_j^{\ul{x}})}$
and $p_{[0,l)}$.
Then $\mathcal{M}_T^p([0,1])$ is dense in
$\mathcal{M}_T^e([0,1])$.
\end{thm}


Hereafter,
let $T\colon [0,1]\to [0,1]$ be a generalized
$(\alpha,\beta)$-transformation and $k$ be a smallest integer
not less than $\alpha+\beta$.
To simplifies the notation, we set
$A_n:=\sigma^n([a_{[0,n]}])$ and $B_n:=\sigma^n([b_{[0,n]}])$ for each $n\ge 0$
and set
$R_m:=R_m^{\ul{a}}$, $S_m:=R_m^{\ul{b}}$,
$r_{m+1}:=r_{m+1}^{\ul{a}}=R_{m+1}-R_m$ and
$s_{m+1}:=r_{m+1}^{\ul{b}}=S_{m+1}-S_m$ for each $m\ge 0$.
We also denote

\vspace{11pt}\noindent
$\mathcal{A}_1:=\{m\in\mathbb{N}:
\sigma(f_m(\ul{a}))=\ul{a}\}$,\hspace{0.5cm}
$\mathcal{A}_2:=\{m\in\mathbb{N}:
\sigma(f_m(\ul{a}))=\ul{b}\}$,\vspace{0.2cm}\\
$\mathcal{B}_1:=\{m\in\mathbb{N}:
\sigma(f_m(\ul{b}))=\ul{a}\}$ and
$\mathcal{B}_2:=\{m\in\mathbb{N}:
\sigma(f_m(\ul{b}))=\ul{b}\}$.

\vspace{11pt}\noindent
Note that $\sigma [i]=\Sigma_T^+$ for each
$2\le i\le k-1$ and $\sigma(\ul{x})\in\{\ul{a},\ul{b}\}$
for any $\ul{x}\in\mathcal{CR}$.
These together with Proposition \ref{diagram-3}
and Theorem \ref{HR-thm}
imply the following:

\begin{prop}
\label{diagram-2}
Let $(\mathcal{D},\rightarrow)$ be the Markov diagram
of $T$.
\vspace{0.2cm}\\
(1) $\mathcal{D}=\{A_n, B_n: n\geq0\}\cup\{[2],\ldots, [k-1]\}$.
\vspace{0.2cm}\\
(2) (i) For any $n\ge 0$, $A_n\rightarrow A_{n+1}$ and $B_n\rightarrow B_{n+1}$.
\vspace{0.2cm}\\
(ii) For $m\geq1$,
\begin{itemize}
\item
$A_{R_m-1}\rightarrow A_{r_m-1}$
and $a_{(R_{m-1},R_m)}=a_{[0,r_m-1)}$ if $m\in\mathcal{A}_1$,
\item
$A_{R_m-1}\rightarrow B_{r_m-1}$
and $a_{(R_{m-1},R_m)}=b_{[0,r_m-1)}$
if $m\in\mathcal{A}_2$.
\end{itemize}
(iii) For $m\geq1$,
\begin{itemize}
\item
$B_{S_m-1}\rightarrow A_{s_m-1}$
and $b_{(S_{m-1},S_m)}=a_{[0,r_m-1)}$ if $m\in\mathcal{B}_1$,
\item
$B_{S_m-1}\rightarrow B_{s_m-1}$
and $b_{(S_{m-1},S_m)}=b_{[0,s_m-1)}$
if $m\in\mathcal{B}_2$.
\end{itemize}

\noindent
(3) There exist two maps $P\colon\mathcal{A}_2\to\mathbb{N}_0$
and $Q\colon\mathcal{B}_1\to\mathbb{N}_0$
such that for any $m\ge 1$, $r_m-1=S_{p(m)}$.
and $s_m-1=R_{Q(m)}$.
\end{prop}

\begin{thm}
\label{HR-thm2}
Suppose that there are integers $N_0$ and
$N_1$ satisfying the following conditions:
\begin{itemize}
    \item 
For any $j\in\mathbb{N}$ with $R_j>N_0$,
there is an integer $1\le m<j$,
a periodic point $\ul{p}$ with period $l$
and $u\in\mathcal{L}(\Sigma_T^+)$ with
$|u|\ge R_j-R_m-N_1$
such that $u$ is a subword of both
$a_{[R_m,R_j)}$ and $p_{[0,l)}$.

\item
For any $j\in\mathbb{N}$ with $S_j>N_0$,
there is an integer $1\le m<j$,
a periodic point $\ul{p}$ with period $l$
and $u\in\mathcal{L}(\Sigma_T^+)$ with
$|u|\ge S_j-S_m-N_1$
such that $u$ is a subword of both
$b_{[R_m,R_j)}$ and $p_{[0,l)}$.
\end{itemize}
Then $\mathcal{M}^p_T([0,1])$ is dense
in $\mathcal{M}_T^e([0,1])$.
\end{thm}


\section{Proofs}

\subsection{Proof of Proposition A}

In this subsection, we give a proof of Proposition A
with the assumption that $T$ is transitive.
As we mentioned in \S1, this theorem appears as
\cite[Theorem A]{CY} with the stronger assumption that
$\mathcal{D}$ is irreducible.
In \cite{CY},
the hypothesis of the irrecducibility for $\mathcal{D}$
is used only in \cite[Proposition 3.1]{CY}.
The other part of the proof of \cite[Theorem A]{CY}
can be shown similarly
by using the transitivity of $T$ instead of the irreducibility.
Hence to prove Proposition A, it is sufficient to show
the following proposition, which is analogous
to \cite[Proposition 3.1]{CY}.

%

\begin{prop}
\label{e-dense}
Let $\mathcal{C}\subset \mathcal{D}$ be as in
Lemma \ref{tran}.
For any $\epsilon>0$,
any $\mu\in\mathcal{M}_{\sigma}(\Sigma_T^+)$,
and any neighborhood $\mathcal{U}\subset \mathcal{M}(\Sigma_T^+)$
of $\mu$,
there exist a finite set $\mathcal{F}\subset \mathcal{C}$ and
$\rho\in\mathcal{M}_{\sigma}^e(\Psi(\Sigma^+_{\mathcal{F}}))$ such that
$\rho\in \mathcal{U}$ and $\textcolor{black}{h_{\sigma}(\rho)}\ge \textcolor{black}{h_{\sigma}(\mu)}-2\epsilon$.
\begin{proof}
Without loss of generality, we may assume
that $h_{\sigma}(\mu)-\epsilon>0$, otherwise the conclusion is yield by the assumption that $\mathcal{M}_{\sigma}^p(\Sigma_T^+)$ is dense
in $\mathcal{M}_{\sigma}^e(\Sigma_T^+)$.
By Ergodic Decomposition Theorem and
the affinity of the entropy map,
there exists a finite convex combination of ergodic measures
$\nu:=\sum_{i=1}^pa_i\nu_i$ such that
$h_{\sigma}(\nu)\in\mathcal{U}$ and
$h_{\sigma}(\nu)\ge h_{\sigma}(\mu)-\epsilon$.
Again by density of $\mathcal{M}^p_{\sigma}(\Sigma_T^+)$
in $\mathcal{M}_{\sigma}^e(\Sigma_T^+)$,
we may assume $\nu_i\in\mathcal{M}_{\sigma}^p(\Sigma_T^+)$
whenever $h_{\sigma}(\nu_i)=0$.
Then we need the following lemma.
\begin{lem}
\label{liftable}
For each $1\le i\le p$, there exists
$\overline{\nu_i}\in\mathcal{M}_{\sigma}^e(\Sigma_{\mathcal{C}}^+)$ such that
$\nu_i=\overline{\nu_i}\circ\Psi^{-1}$ holds.
In particular, $\overline{\nu}:=\sum_{i=1}^pa_i
\overline{\nu_i}$ satisfies
$\overline{\nu}\in\mathcal{M}_{\sigma}
(\Sigma_{\mathcal{C}}^+)$
and $\nu=\overline{\nu}\circ \Psi^{-1}$.
\begin{proof}
We divide the proof into two cases.
\vspace{0.2cm}

{\bf(Case 1)} $h_{\sigma}(\nu_i)=0$.
In this case, $\nu_i$ is a periodic measure.
Take a periodic point $\ul{x}\in \Sigma_T^+$ in the support of $\nu_i$.
Then it follows from \cite[Theorem 8]{H3} and
$\Psi(\Sigma^+_{\mathcal{C}})=\Sigma_T^+$ that
there are finite vertices $C_0,\ldots,C_{n-1}\in\mathcal{C}$
such that $(C_{[0,n)})^{\infty}\in\Sigma_{\mathcal{C}}^+$
and $\Psi((C_{[0,n)})^{\infty})=\ul{x}$.
Hence if we set $\overline{\nu_i}:=
\delta_n^{\sigma}((C_{[0,n)})^{\infty}),$
then we have $\overline{\nu_i}\in\mathcal{M}_{\sigma}^e
(\Sigma_{\mathcal{C}}^+)$
and $\nu_i=\overline{\nu_i}\circ \Psi^{-1}$.

{\bf (Case 2)} $h_{\sigma}(\nu_i)>0$.
Since $\nu_i\in\mathcal{M}_{\sigma}^e(\Sigma_T^+)$
and $h_{\sigma}(\nu_i)>0$,
it follows from Lemma \ref{positive} that
there exists $\overline{\nu_i}\in\mathcal{M}_{\sigma}^e(\Sigma_{\mathcal{D}}^+)$
such that
$\overline{\nu_i}:=\nu_i\circ\Psi^{-1}$ holds.
Since $\overline{\nu_i}$ is ergodic, there exists
an irreducible subset $\mathcal{C}'\subset \mathcal{D}$
such that $\overline\nu_i(\Sigma_{\mathcal{C}'}^+)=1$.
Assume that
$\mathcal{C}'\not=\mathcal{C}$.
Since $\Sigma^+_{\mathcal{C}'}\subset \Psi^{-1}(\Psi(\Sigma^+_{\mathcal{C}'}))$, we have
$\nu_i(\Psi(\Sigma^+_{\mathcal{C}'}))=
\overline{\nu_i}\circ \Psi^{-1}
(\Psi(\Sigma^+_{\mathcal{C}'}))\ge \overline{\nu_i}
(\Sigma^+_{\mathcal{C}'})=1$.
This implies that $\nu_i(\Psi(\Sigma^+_{C})\cap
\Psi(\Sigma^+_{\mathcal{C}'}))=\nu_i(\Sigma_T^+\cap \Psi(\Sigma_{\mathcal{C}'}^+))=1$.
On the other hand, by \cite[Theorem 1 (ii)]{H4},
$\Psi(\Sigma^+_{C})\cap
\Psi(\Sigma^+_{\mathcal{C}'})$ is either empty or finite
(see also \cite[Page 385]{H3}).
Hence by the ergodicity of $\nu_i$, we have $h_{\sigma}(\nu_i)=0$,
which is a contradiction.
\end{proof}
\end{lem}
Note that $\overline{\nu}$ is an invariant measure on
a transitive countable Markov shift $\Sigma^+_{\mathcal{C}}$.
Hence by the continuity of $\Psi$
and \cite[Main Theorem]{T}, we can find a finite set
$\mathcal{F}\subset \mathcal{C}$ and
an ergodic measure
$\overline{\rho}$ on $\Sigma^+_{\mathcal{F}}$ such that
$h_{\sigma}(\overline{\rho})\ge h_{\sigma}(\overline{\nu})
-\epsilon$ and
$\rho:=\overline{\rho}\circ \Psi^{-1}\in\mathcal{U}$.
Since $\Psi\colon \Sigma^+_{\mathcal{C}}\to\Sigma_T^+$
is countable to one, by \cite[Proposition 2.8]{Bu3},
we have $h_{\sigma}(\nu)=h_{\sigma}(\overline{\nu})$ and
$h_{\sigma}(\rho)=h_{\sigma}(\overline{\rho})$,
which prove the proposition.
\end{proof}
\end{prop}

\subsection{Proof of Theorem \ref{diagram}}
The aim of this subsection is to give a proof of
Theorem \ref{diagram}.
Let $n_0$ be as in Lemma \ref{tran} and
$m_0$ be a smallest integer such that
$R_{m_0}, S_{m_0}\ge n_0$ holds.
We set $N_0:=\max\{R_{m_0},\ S_{m_0}\}$ and
let $n_1$ be a smallest number such that
for any $C\in\mathcal{C}\cap\mathcal{D}_{n_0}$
and $D\in\mathcal{C}\cap\mathcal{D}_{N_0}$,
there exist finite vertices $C_0,\ldots,C_n$
with $n\le n_1$ such that
$C_{[0,n]}\in\mathcal{L}(\Sigma_{\mathcal{C}}^+)$,
$C_0=C$ and $C_n=D$.
We note that $n_1<\infty$ since $\mathcal{C}$ is
irreducible. We set $N_1:=N_0+n_1$.
By Theorem \ref{HR-thm2}, to prove Theorem \ref{diagram},
it is sufficient to show the following:
\begin{enumerate}
\item 
For any $j\in\mathbb{N}$ with $R_j>N_0$,
there is an integer $1\le m<j$,
a periodic point $\ul{p}$ with period $l$
and $u\in\mathcal{L}(\Sigma_T^+)$ with
$|u|\ge R_j-R_m-N_1$
such that $u$ is a subword of both
$a_{[R_m,R_j)}$ and $p_{[0,l)}$.

\item
For any $j\in\mathbb{N}$ with $S_j>N_0$,
there is an integer $1\le m<j$,
a periodic point $\ul{p}$ with period $l$
and $u\in\mathcal{L}(\Sigma_T^+)$ with
$|u|\ge S_j-S_m-N_1$
such that $u$ is a subword of both
$b_{[R_m,R_j)}$ and $p_{[0,l)}$.
\end{enumerate}
We only prove the item (I) because (II)
can be shown in a similar manner.
Take any $j\in\mathbb{N}$ with $R_j>N_0$.
In what follows we will decompose
the set $\mathbb{N}$ into six sets.
Let $\mathcal{A}_1$, $\mathcal{A}_2$, $\mathcal{B}_1$
and $\mathcal{B}_2$ be as in \S2.3 and
set

\vspace{11pt}\noindent
$\mathcal{A}_3:=\{m\in\mathbb{N}:\mathcal{S}(A_{R_m-1})
\cap\mathcal{D}_0\not=\emptyset\}$
and\vspace{0.2cm}\\
$\mathcal{B}_3:=\{m\in\mathbb{N}:\mathcal{S}(B_{S_m-1})\cap
\mathcal{D}_0\not=\emptyset\}$.

\vspace{11pt}\noindent
First, it is clear that
$\mathbb{N}=\mathcal{A}_1\cup\mathcal{A}_2
=\mathcal{A}_1\cup (\mathcal{A}_2\setminus\mathcal{A}_3)
\cup \mathcal{A}_3$.
We set

\vspace{11pt}\noindent
$\mathcal{A}_2^{(1)}:=\{m\in\mathcal{A}_2:
m+1\in\mathcal{A}_1\},$\ \ \ \ \ \ \ \ \ \ \ 
$\mathcal{A}_2^{(2)}:=\{m\in\mathcal{A}_2:
P(m)\in\mathcal{B}_3\},$\vspace{0.2cm}\\
$\mathcal{A}_2^{(3)}:=\{m\in\mathcal{A}_2:
P(m)+1\in\mathcal{B}_2\},$ and
$\mathcal{A}_4:=\mathcal{A}_2\setminus
(\mathcal{A}_3\cup\bigcup_{j=1}^3\mathcal{A}_2^{(j)})$,

\vspace{11pt}\noindent
noting that the map $P$ is
defined on $\mathcal{A}_2$
by Proposition \ref{diagram-2} (3).
\\
Then we have $\mathcal{A}_2\setminus\mathcal{A}_3
\subset\bigcup_{j=1}^3\mathcal{A}_2^{(j)}\cup \mathcal{A}_4$, which
implies that
\[\mathbb{N}=\mathcal{A}_1\cup
\bigcup_{j=1}^3\mathcal{A}_2^{(j)}\cup\mathcal{A}_3\cup
\mathcal{A}_4.\]
Hence we can divide the proof into six cases:

\vspace{11pt}
\noindent
{\bf (Case A) $j\in\mathcal{A}_1$.\hspace{0.65cm}
(Case B) $j\in\mathcal{A}_2^{(1)}$.\hspace{0.65cm}
(Case C) $j\in\mathcal{A}_2^{(2)}$.\vspace{0.2cm}\\
(Case D) $j\in\mathcal{A}_2^{(3)}$.\hspace{0.4cm}
(Case E) $j\in\mathcal{A}_3$.\hspace{0.92cm}
(Case F) $j\in\mathcal{A}_4.$}

\vspace{11pt}\noindent
We note that the proofs of {\bf (Case A)}
and {\bf (Case B)} are
similar to those of {\bf (Case 1)}
and {\bf (Case 2)} in \cite{HR} respectively.

\vspace{11pt}
{\bf (Case C) $j\in\mathcal{A}_2^{(2)}$.}
Take $D\in\mathcal{S}(B_{S_{P(j)}-1})\cap\mathcal{D}_0$.
By the definition of $n_1$, we can find
an integer $n\le n_1$ and finite vertices
$C_0,\ldots,C_n\in\mathcal{C}$ such that
$C_{[0,n]}\in\mathcal{L}(\Sigma^+_{\mathcal{C}})$,
$C_0=D$ and $C_n=B_{S_{m_0}}$.
Now, we set
$\ul{p}:=\Psi((B_{[S_{m_0},S_{P(j)})}C_{[0,n)})^{\infty})$.
Then $\ul{p}$ is a periodic point with period
$l:=S_{P(j)}-S_{m_0}+n$.
Since $A_{R_j-1}\rightarrow B_{S_{P(j)}}$,
we have $a_{(R_{j-1},R_j)}=b_{[0,S_{P(j)})}$,
which implies (I).

\vspace{11pt}
{\bf (Case D) $j\in\mathcal{A}_2^{(3)}$.}
Since $P(j)+1\in\mathcal{B}_2$,
we can find $u\le P(j)$ such that $B_{S_{P(j)+1}-1}\rightarrow B_u$.
We set $\ul{p}:=\Psi((B_{(S_{P(j)},S_{P(j)+1})}
B_{[u,S_{P(j)}]})^{\infty})$.
Then $\ul{p}$ is a periodic point
with period $l:=S_{P(j)+1}-u$.
Since $A_{R_j-1}\rightarrow B_{S_{P(j)}}$ and
$B_{S_{P(j)+1}-1}\rightarrow B_u$, we have
$a_{(R_{j-1},R_j)}=b_{[0,S_{P(j)})}
=b_{[0,u)}b_{[u,S_{P(j)})}=b_{(S_{P(j)},S_{P(j)+1})}
b_{[u,S_{P(j)}]}$. This implies (I).

\vspace{11pt}
{\bf (Case E) $j\in\mathcal{A}_3$.}
Take $D\in \mathcal{S}(A_{R_j-1})\cap\mathcal{D}_0$. By the definition of $n_1$, we can find
an integer $n\le n_1$ and finite vertices
$C_0,\ldots,C_n\in\mathcal{C}$ such that
$C_{[0,n]}\in\mathcal{L}(\Sigma^+_{\mathcal{C}})$,
$C_0=D$ and $C_n=A_{R_{m_0}}$.
Hence if we set $\ul{p}:=\Psi((A_{[R_{m_0},A_{R_j})}C_{[0,n)})^{\infty})$, then (I) holds.

\vspace{11pt}
{\bf (Case F) $j\in\mathcal{A}_4$.}
In this case,
$j, j+1\in\mathcal{A}_2$, $P(j)+1\in\mathcal{B}_1$, $j\not\in\mathcal{A}_3$
and $P(j)\not\in\mathcal{B}_3$ hold
by the definition of $\mathcal{A}_4$.
We prove the following proposition,
which plays a fundamental role to prove this case.

\begin{prop}
\label{no long edges}
We have
\[    S_{P(j)}^\infty \succ \{r_{j+i}-1\}_{i=1}^{\infty}
        \quad \mbox{or}\quad
        R_j^\infty \succ \{s_{P(j)+i}-1\}_{i=1}^{\infty}.
\]
Here $\succ$ denotes the lexicographical order.

\begin{proof}
For the notational simplicity, set
$r^{(m)}:=\{r_{m+i}-1\}_{i=1}^{\infty}$ and
$s^{(m)}:=\{s_{m+i}-1\}_{i=1}^{\infty}$ for $m\ge 1$.
Showing by contradiction, we assume
\begin{equation}
    S_{P(j)}^\infty \preceq r^{(j)}\quad \mbox{and}\quad R_j^\infty \preceq s^{P(j)}.
    \label{contradiction}
\end{equation}
For $q\geq R_j+1$ set $\mathcal{E}_{q,1}:=\{m: R_j\leq R_m <q\}$ and $\mathcal{E}_{q,2}:=\{m: S_{P(j)}\leq R_m <q\}$.
Denote also $\ul{x}:=f_{j+1}(\ul{a})$ and
$\ul{y}:={\rm adj}(\ul{x})=f_{P(j)+1}(\ul{b})$.
Note that $f_m(\ul{a})=x$ implies $m\in\mathcal{A}_2$.
Similarly $f_m(\ul{b})=\ul{y}$ implies
 $m\in\mathcal{B}_1$.
To prove the proposition, it is sufficient to show
the following for all $q\ge R_j+1$:
\begin{itemize}
\item For any $m\in\mathcal{E}_{q,1}$, we have
\begin{equation}
\label{A}    
m\not\in\mathcal{A}_3,\ f_{m+1}(\ul{a})=\ul{x}\text{ and }S_{P(j)}^{\infty}\preceq r^{(m)}.
\end{equation}

\item
For any $m\in\mathcal{E}_{q,2}$, we have
\begin{equation}
\label{B}
m\not\in\mathcal{B}_3,\ f_{m+1}(\ul{b})=\ul{y}\text{ and }R_j^{\infty}\preceq s^{(m)}. 
\end{equation}
\end{itemize}
Indeed, letting $\mathcal{C}':=\{A_m :m\ge R_j\}\cup\{B_m: m\ge S_{P(j)}\}$, 
we have
$\bigcup_{C\in\mathcal{C}'}\mathcal{S}(C)\cap (\mathcal{C}\setminus\mathcal{C}')=\emptyset$
by \eqref{A} and \eqref{B}.
However, this is a contradiction
since $\mathcal{C}$ is irreducible and $A_{n_0}\in
\mathcal{C}\setminus\mathcal{C}'$.

We prove \eqref{A} and \eqref{B} by induction with regard to $q\ge R_j+1$.
Let $q=R_j+1$.
Since $R_{j+1}\geq R_j+1$, $j+1\notin \mathcal{E}_{R_j+1,1}$.
Similarly, $P(j)+1\notin \mathcal{E}_{R_j+1,2}$ because $S_{P(j)+1}>s_{P(j)+1}\geq R_j$.
Hence we have $\mathcal{E}_{R_{j}+1,1}=\{j\}$ and $\mathcal{E}_{R_{j}+1,2}=\{P(j)\}$.
Hence \eqref{A} and \eqref{B} immediately follows by
\eqref{contradiction}.

\vspace{11pt}
Let $q>R_j+1$ and assume the inductive hypothesis
for $q-1$.
Take $m\in \mathcal{E}_{q,1}\setminus \mathcal{E}_{q-1,1}$
(if no such $m$ exists, \eqref{A} automatically holds for all $m\in \mathcal{E}_{q,1}$).
We will show that \eqref{A} holds for $m$.
Since $m\in\mathcal{E}_{q,1}$, we have
$R_{m-1}<R_m<q$ and
hence $m-1\in \mathcal{E}_{q-1,1}$.
Therefore, by the inductive hypothesis, we have
$m-1\notin\mathcal{A}_3$, $f_m(\ul{a})=\ul{x}$ and
$S_{P(j)}^\infty \preceq r^{(m-1)}$.
%
\begin{figure}
    \centering
    \includegraphics[width=200pt]{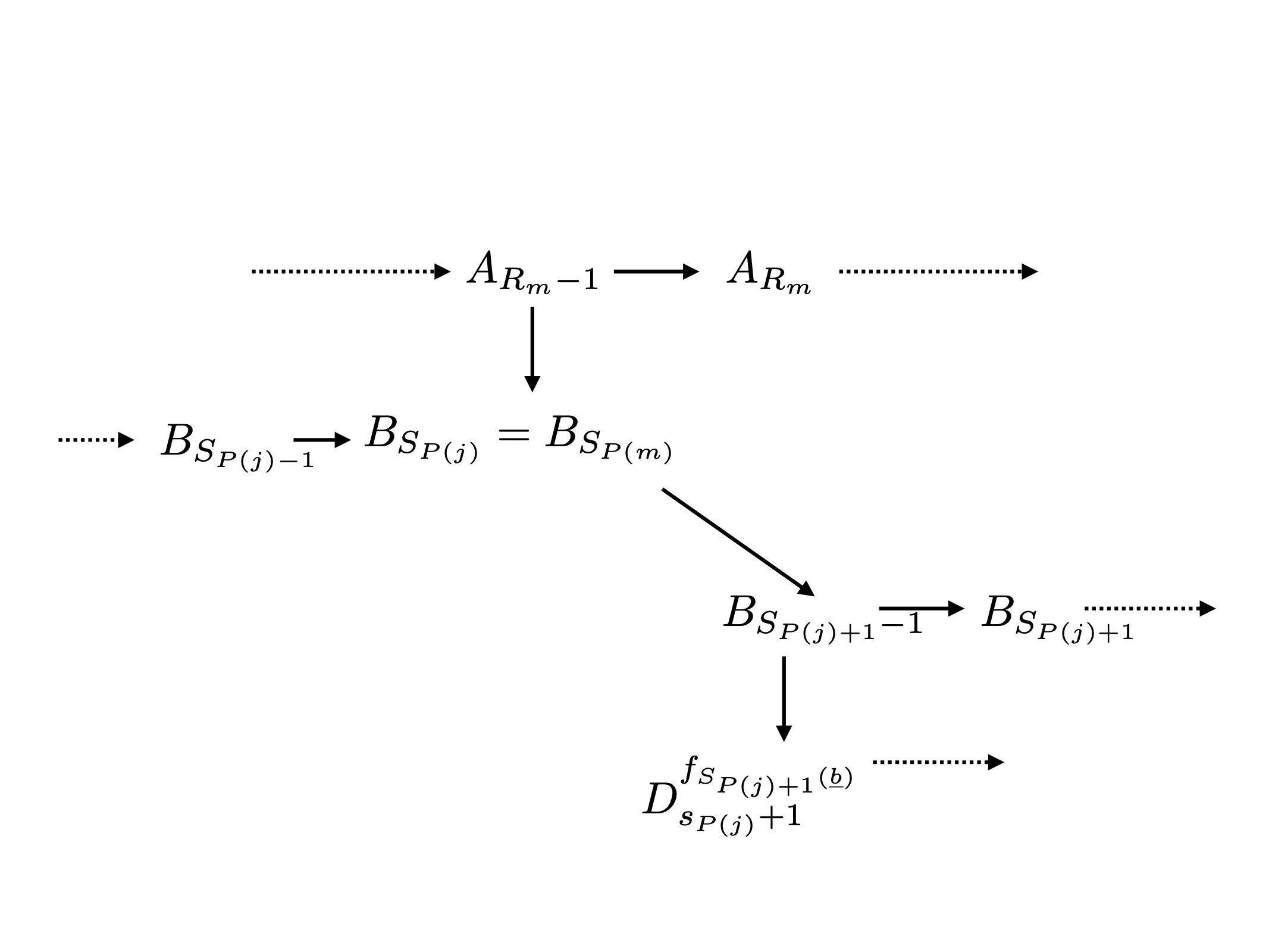}
    \caption{Sketch of the Markov diagram in Case 1.}
    \label{graph_case1}
\end{figure}

\vspace{11pt}
{\bf (Case 1)} $ S_{P(j)}=r_m-1$ (for the situation of the Markov diagram, see Figure \ref{graph_case1}).
In this case, we have $S_{P(j)}=r_m-1=S_{P(m)}$.
This implies $A_{R_m-1}\subset B_{S_{P(j)-1}}$ and $\#\mathcal{S}(A_{R_m-1})\leq \#\mathcal{S}(B_{S_{P(j)-1}})=2$.
Hence we have $m\notin \mathcal{A}_3$.
Moreover, 
$f_{m+1}(\ul{a})=\ul{x}$ is followed by 
$\mathcal{S}(A_{R_m-1})=\{A_{R_m}, B_{S_{P(j)}}\}$ and $f_{S_{P(j)+1}}(\ul{b})=\ul{y}$.
The last condition $S_{P(j)}^\infty \preceq
r^{(m)}$
follows from $S_{P(j)}=r_m-1$ and
$S_{P(j)}^\infty \preceq r^{(m-1)}$.

\vspace{11pt}
{\bf (Case 2)} $S_{P(j)}<r_m-1$ (For the situation of the Markov diagram, see Figure \ref{graph_case2}).
In this case, we have
$S_{P(j)}<r_m-1=S_{P(m)}<R_m=q-1$, which implies that
$P(m)-1,\ P(m)\in \mathcal{E}_{q-1,2}$.
Hence by the inductive hypothesis (\ref{B}) for $q-1$,
we have
$P(m)-1, P(m)\notin \mathcal{B}_3$, $P(m), P(m)+1 \in \mathcal{B}_1$ and $R_j^\infty \preceq s^{(P(m)-1)},\ s^{(P(m))}$.
Since $A_{R_m-1}\rightarrow B_{S_{P(m)}}$, 
we have $\# \mathcal{S}(A_{R_m-1})\leq \#\mathcal{S}(B_{S_{P(m)-1}})=2$,
which implies $m\notin\mathcal{A}_3$.
%
Moreover,
$\mathcal{S}(A_{R_m-1})=\{A_{R_m}, B_{S_{P(m)}}\}$ implies
\[f_{m+1}(\ul{a})={\rm adj}(f_{S_{P(m)+1}}(\ul{b}))={\rm adj}(\ul{y})=\ul{x}.\]

Now we consider the last condition $S_{P(j)}^\infty\preceq 
r^{(m)}$.
Since 
$Q(P(m))=s_{P(m)}-1\geq R_j$ and $Q(P(m))< s_{P(m)}<S_{P(m)}<R_m=q-1,$
\eqref{A} holds for $Q(P(m))$.
If $r^{(Q(P(m)))} = r^{(m)}$, \eqref{A} yields
$S_{P(j)}^\infty\preceq r^{(Q(P(m)))}= r^{(m)}$.
Otherwise let $\ell=\inf\{i\geq 1: r_{Q(P(m))+i}\neq r_{m+i}\}$. 
It is easy to see $S_{P(j)}^\infty\preceq r^{(m)}$ is yield by the following:
For $i=0, 1, \ldots, \ell$, we have
\begin{equation}
    A_{R_{m+i}} \subset A_{R_{Q(P(m))+i}} 
    \ \mbox{and}\ 
    \#\tau(A_{R_{m+i}-1})=\#\tau(A_{R_{Q(P(m))+i}-1} )=2 \label{two ss}
\end{equation}
In particular, \eqref{two ss} implies $r_{Q(P(m))+\ell+1} < r_{m+\ell+1}$ and
$r^{(Q(P(m))+1)}\preceq r^{(m)}$.
We close our proof of (\ref{A}) by showing \eqref{two ss}.
Since $A_{R_m-1}\rightarrow  B_{S_{P(m)}}$ and $B_{S_{P(m)}-1}\rightarrow A_{R_{Q(P(m))}}$, we have
\[
    A_{R_m-1}\subset B_{S_{P(m)}-1}\subset  A_{R_{Q(P(m))}-1}.
\]
Moreover, 
$f_{m+1}(a)=f_{Q(P(m))+1}=\ul{x}$ implies $A_{R_m}, A_{R_{Q(P(m))}} \subset [x_0]$ and 
$A_{R_m}\subset A_{R_{Q(P(m))}}$.
Since \eqref{A} holds for $Q(P(m))$, we have
$\#\mathcal{S}(A_{R_m-1})=\#\mathcal{S}(A_{R_{Q(P(m))}-1} )=2$.
%
%
\begin{figure}
    \begin{center}
    \includegraphics[width=200pt]{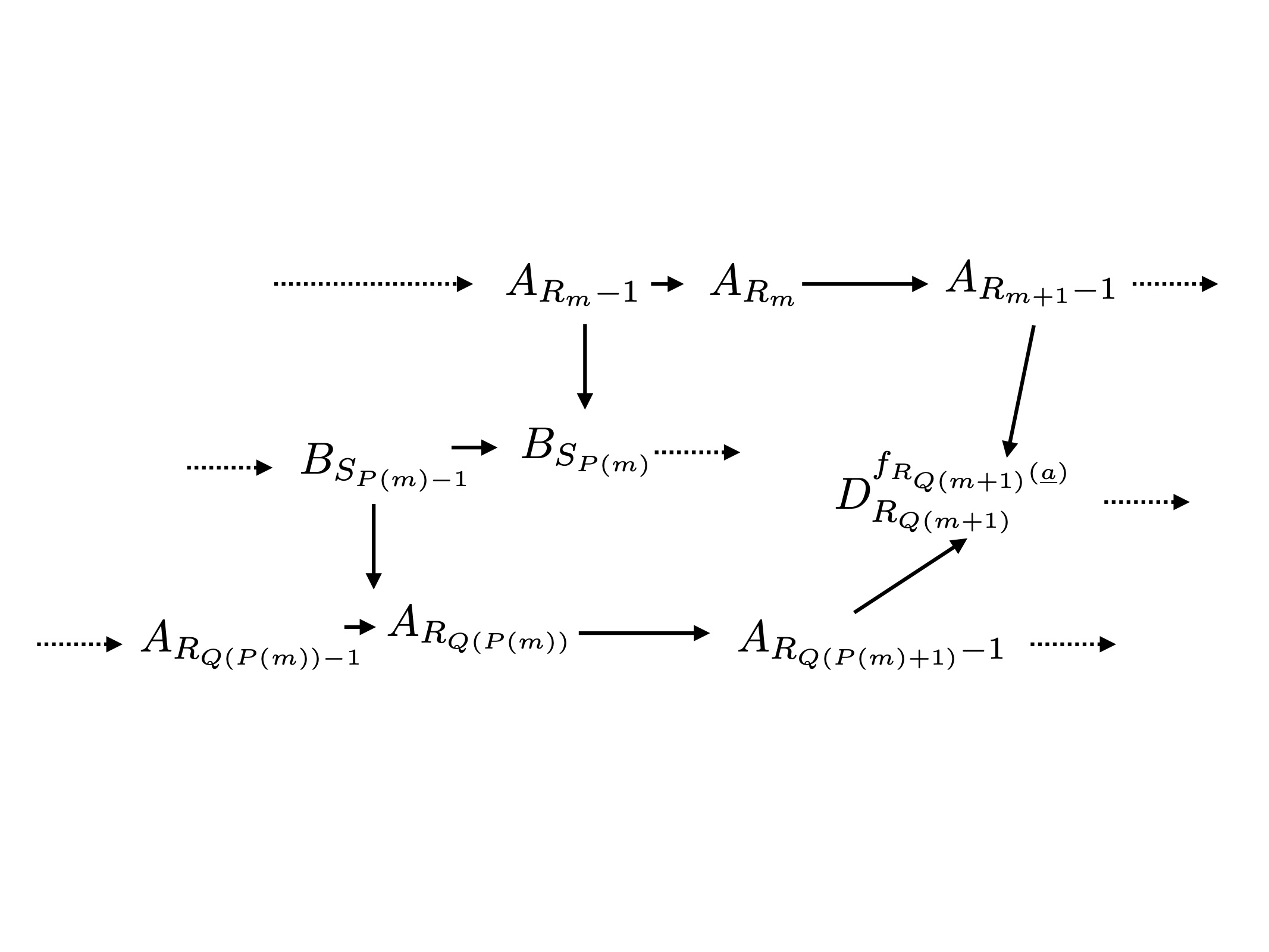}
    \caption{Sketch of the Markov diagram in Case 2.}
    \end{center}
    \label{graph_case2}
\end{figure}

Let $0<i\leq \ell$ and
assume \eqref{two ss} holds for all $0\leq i'<i$.
Since $A_{R_m+i-1}\subset A_{R_{Q(P(m))}+i-1}$ and $r_{Q(P(m))+i-1} = r_{m+i-1}$,
we have
$A_{R_{m+i}-1}\subset A_{R_{Q(P(m))+i}-1}$.
Since $A_{R_{Q(P(m))}}<S_{P((m)}<R_m$ yields 
\begin{align*}
    R_{Q(P(m))+i}&=R_{Q(P(m))}+\sum_{i'=1}^i r_{Q(P(m))+i'}\\
    &=R_{Q(P(m))}+\sum_{i'=1}^i r_{m+i'}\\
    &<R_m+R_{m+i}-R_m=R_{m+i},
\end{align*}
we have
$\#\mathcal{S}(A_{R_{m+i}-1})=\#\mathcal{S}(A_{R_{Q(P(m))+i}-1})=2$.
Combining the definition of $\ell$ and the \eqref{two ss} for $\ell$, 
we have 
$r_{Q(P(m))+\ell+1} < r_{m+\ell+1}$.

We can prove (\ref{B}) for all $l\in\mathcal{E}_{q,2}$
in a similar manner,
which proves the proposition.
\end{proof}
\end{prop}

We continue to the proof of Theorem \ref{diagram}.
By Proposition \ref{no long edges}, we can divide
(Case F) into the following three cases:

\vspace{11pt}\noindent
{\bf (Case F1)} $R_j^{\infty}\succ \{s_{P(j)+i}-1\}_{i=1}^{\infty}$.
\vspace{0.2cm}\\
{\bf (Case F2)} $S_{P(j)}^\infty \succ \{r_{j+i}-1\}_{i=1}^\infty$ and $P(j)-1\in \mathcal{B}_1$.
\vspace{0.2cm}\\
{\bf (Case F3)} $S_{P(j)}^\infty \succ \{r_{j+i}-1\}_{i=1}^\infty$ and $P(j)-1\in\mathcal{B}_2$.
\vspace{11pt}

The proofs of {\bf (Case F1)}, {\bf (Case F2)} and
{\bf (Case F3)} are similar to those of
{\bf (Case 3)}, {\bf (Case 4)} and {\bf (Case 5)}
in \cite{HR}, respectively.
Theorem \ref{diagram} is proved.

\begin{remark}\label{novelty}
In our proof of Theorem \ref{diagram}, we improve the method in Hofbauer and Raith's work \cite{HR} to apply generalized $(\alpha, \beta)$-transformations.
For a piecewise monotonic map with two monotonic pieces,
its critical set $\mathcal{CR}$ consists of two points which are clearly adjacent.
Moreover, its Markov diagram has no vertices at which more than three edges start
and a word represented by a path without branching vertices coincides with the word of either of the points in $\mathcal{CR}$.
On the other hand, for a generalized $(\alpha, \beta)$-transformation, 
there may exist a vertex at which more than three edges start in its Markov diagram.
Moreover it is difficult to check which critical point represent a word defined on a path without branching vertices.
Therefore, in the proof of Proposition \ref{no long edges} we need to be careful to check that every vertex has two or less edges and a pair of adjacent critical points, $\ul{x}=f_j(\ul{a})$ and $\ul{y}={\rm adj}(\ul{x})=f_{P(j+1)}(\ul{b})$, appears in each inductive step. 
\end{remark}

\vspace*{3mm}

\noindent
\textbf{Acknowledgement.}~ 
The authors would like to thank
Hajime Kaneko for suggesting this problem.
The first author was partially supported by JSPS KAKENHI Grant Number 21K13816 and
the second author was partially supported by
JSPS KAKENHI Grant Number 21K03321.

\end{document}